\documentclass{siamart190516}

\usepackage{lipsum}
\usepackage{amsfonts}
\usepackage{graphicx}
\usepackage{epstopdf}
\usepackage{algorithmic}
\usepackage{multirow}
\usepackage{booktabs}

\usepackage{tikz}
\usetikzlibrary{arrows}
\usetikzlibrary{decorations.pathreplacing}
\usetikzlibrary{shapes}

\ifpdf%
  \DeclareGraphicsExtensions{.eps,.pdf,.png,.jpg}
\else
  \DeclareGraphicsExtensions{.eps}
\fi
\usepackage{amsopn}

\usepackage{booktabs}

\renewcommand{\vec}[1]{\boldsymbol #1}
\newcommand{\vv}{\vec{v}}
\newcommand{\vu}{\vec{u}}
\newcommand{\ve}{\vec{e}}
\newcommand{\vr}{\vec{r}}
\newcommand{\vf}{\vec{f}}
\newcommand{\vg}{\vec{g}}
\newcommand{\vtau}{\vec{\tau}}
\newcommand{\fg}{\Omega_f}
\newcommand{\cg}{\Omega_c}

\newcommand{\TheTitle}{%
  Toward Parallel in Time for Chaotic Dynamical Systems
}

\newcommand{\TheShortTitle}{%
  \TheTitle
}

\newcommand{\TheName}{%
  David A. Vargas
}

\newcommand{\TheAddress}{%
  University of New Mexico,
  (\email{dvargas2@unm.edu}).
}

\newcommand{\TheFunding}{%
This work performed under the auspices of the US Department of
Energy by Lawrence Livermore National Laboratory under Contract
DE-AC52-07NA27344 (LLNL-CONF-830992).
}

\newcommand{\TheCollaborators}{%
Robert D. Falgout,
Stefanie G{\"u}nther, and
Jacob B. Schroder
}

\author{\TheName\thanks{\TheAddress}}
\title{{\TheTitle}\thanks{\TheFunding}}
\headers{\TheShortTitle}{\TheName}
\ifpdf%
\hypersetup{%
  pdftitle={\TheTitle},
  pdfauthor={\TheName}
}
\fi

\begin{document}

\maketitle

\begin{center}
  In collaboration with:
  {\TheCollaborators}
\end{center}
\vspace{1cm}

\begin{abstract}
As CPU clock speeds have stagnated, and high performance computers continue to have ever higher core counts, increased parallelism is needed to take advantage of these new architectures. Traditional serial time-marching schemes are a significant bottleneck, as many types of simulations require large numbers of time-steps which must be computed sequentially. Parallel in Time schemes, such as the Multigrid Reduction in Time (MGRIT) method, remedy this by parallelizing across time-steps, and have shown promising results for parabolic problems. However, chaotic problems have proved more difficult, since chaotic initial value problems are inherently ill-conditioned. MGRIT relies on a hierarchy of successively coarser time-grids to iteratively correct the solution on the finest time-grid, but due to the nature of chaotic systems, subtle inaccuracies on the coarser levels can lead to poor coarse-grid corrections. Here we propose a modification to nonlinear FAS multigrid, as well as a novel time-coarsening scheme, which together better capture long term behavior on coarse grids and greatly improve convergence of MGRIT for chaotic initial value problems. We provide supporting numerical results for the Lorenz system model problem.
\end{abstract}


\section{Introduction}\label{sec:intro}
While Parallel in Time (PinT) methods date back over 50 years, interest in these methods has only recently picked up due to the stagnation of CPU clock speeds in the 2000s \cite{gander_review,ben_schroder_review}. For many problems, spatial parallelism can become exhausted, while the time dimension remains largely unparallelized. Thus, we see that PinT schemes, if perfected, have a huge potential for speedup when combined with existing spatial parallel techniques. The reason that the potential of PinT methods has yet to be realized is that the time dimension presents difficulties not seen in the spatial dimension. The most important difference is \emph{causality,} since the solution at a later time depends on the solution at previous times. PinT has already been demonstrated to provide substantial speedups for parabolic problems \cite{parabolic_mgrit}, such as the heat equation, largely due to the fact that the causality of the system is relaxed over time. Parabolic problems have \emph{weak} dependence on initial conditions, since they tend toward a steady state which is largely uncorrelated to the initial data. Hyperbolic problems, such as the wave equation, remain difficult to parallelize in time, since they have \emph{strong} dependence on initial conditions, although some speedup has been demonstrated for such problems in special cases \cite{ben_schroder_review}. However, to our knowledge, no speedup has been achieved for chaotic problems, which exhibit \emph{sensitive} dependence on the initial condition, with the result that the initial value problem is ill-conditioned for chaotic systems. Despite these difficulties, hyperbolic and chaotic systems are very important classes of problems with a wide range of applications across science and engineering.

In order to solve chaotic problems with PinT, we consider multigrid, due to its parallel scalability and potential optimality. However, we note there are many other promising PinT approaches, direct and iterative\cite{gander_review,ben_schroder_review}. We will motivate a modification to the nonlinear FAS coarse grid equation as well as a modified time-coarsening scheme for time-stepping propagators which together greatly improve the convergence of Multigrid Reduction in Time (MGRIT) \cite{MGRIT14} for the chaotic Lorenz system. First we introduce the standard MGRIT algorithm, and study its performance on the Lorenz system. Then we present and motivate the modifications, followed by supporting numerical results.

\subsection{MGRIT}\label{sec:intro_MGRIT}
MGRIT is an iterative multigrid method for solving discrete initial value problems given in the form
\begin{equation}
  \begin{cases}
    \vu_0 = \vf_0                                               \\
    \vu_{i+1} = \Phi(\vu_i) + \vf_{i+1} & i = 0, 1, 2, \dots, n
  \end{cases}
  \label{eqn:ivp}
\end{equation}

where $\Phi$ is a nonlinear time-stepping operator. This system is defined over a finite time-grid with $n + 1$ points, $\fg = \{t_i\}_{i=0}^n$ with time-step size $h = t_{i+1} - t_i$. We will assume, without loss of generality, that $h$ is constant. Let $\vu = \begin{bmatrix} \vu_0,&\vu_1,&\dots,&\vu_n\ \end{bmatrix}^T$ be the state vector, $\vf = \begin{bmatrix} \vf_0, & \vf_1, & \dots, & \vf_n \end{bmatrix}^T$ be a constant forcing term which also encodes the initial condition, then the system of equations \eqref{eqn:ivp} may be written in the form of a block non-linear matrix equation,
\begin{equation}
  A (\vu) = \vf \text{, where } A(\vu) = \begin{bmatrix}
    I      &       &        &        &   \\
    - \Phi & I     &        &        &   \\
           & -\Phi & I      &        &   \\
           &       & \ddots & \ddots &   \\
           &       &        & - \Phi & I
  \end{bmatrix}
  \begin{bmatrix}
    \vu_0 \\ \vu_1 \\ \vu_2 \\ \vdots \\ \vu_n
  \end{bmatrix}
  \label{eqn:systemA}.
\end{equation}
Often, this system comes from a discretization of an ODE of the form $\vu'(t) = \vg(\vu(t))$, in which case $\Phi$ might be e.g. Euler's method. Typically, this system would be solved using forward substitution, which corresponds with time-marching. MGRIT instead applies FAS multigrid to the system \eqref{eqn:systemA}, allowing it to be solved iteratively in parallel. This works by approximately solving \eqref{eqn:systemA} on a hierarchy of coarser time-grids, e.g. $\Omega^{2h}, \Omega^{4h}, \Omega^{8h}, \dots$, and then interpolating error corrections to finer grids, while the finer grids provide further corrections via local relaxation (block Jacobi).

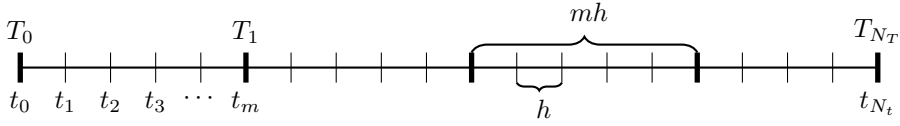
\begin{figure}[t]
  \centering
  \begin{tikzpicture}[xscale = 1.2]
    \draw [thick] (-4.5,0) -- (5.0,0);

    \draw[line width=2pt] (-4.5,-0.2) -- (-4.5,0.2);
    \draw (-4.0,-0.2) -- (-4.0,0.2);
    \draw (-3.5,-0.2) -- (-3.5,0.2);
    \draw (-3.0,-0.2) -- (-3.0,0.2);
    \draw (-2.5,-0.2) -- (-2.5,0.2);
    \draw[line width=2pt] (-2.0,-0.2) -- (-2.0,0.2);
    \draw (-1.5,-0.2) -- (-1.5,0.2);
    \draw (-1.0,-0.2) -- (-1.0,0.2);
    \draw (-0.5,-0.2) -- (-0.5,0.2);
    \draw (0,-0.2) -- (0,0.2);
    \draw[line width=2pt] (0.5,-0.2) -- (0.5,0.2);
    \draw (1.0,-0.2) -- (1.0,0.2);
    \draw (1.5,-0.2) -- (1.5,0.2);
    \draw (2.0,-0.2) -- (2.0,0.2);
    \draw (2.5,-0.2) -- (2.5,0.2);
    \draw[line width=2pt] (3.0,-0.2) -- (3.0,0.2);
    \draw (3.5,-0.2) -- (3.5,0.2);
    \draw (4.0,-0.2) -- (4.0,0.2);
    \draw (4.5,-0.2) -- (4.5,0.2);
    \draw[line width=2pt] (5.0, -0.2) -- (5.0, 0.2);

    \node [below] at (-4.5, -0.2) {$t_0$};
    \node [above] at (-4.5, 0.2) {$T_0$};
    \node [below] at (-4.0, -0.2) {$t_1$};
    \node [below] at (-3.5, -0.2) {$t_2$};
    \node [below] at (-3.0, -0.2) {$t_3$};
    \node [below] at (-2.5, -0.2) {$\cdots$};
    \node [below] at (-2.0, -0.2) {$t_m$};
    \node [above] at (-2.0, 0.2) {$T_1$};
    \node [below] at (5.0, -0.2) {$t_{N_t}$};
    \node [above] at (5.0, 0.2) {$T_{N_T}$};

    \draw [thick,decorate,decoration={brace,amplitude=6pt,raise=0pt}] (0.5, 0.2) -- (3.0, 0.2);
    \node [above] at (1.8, 0.5) {$m h$};
    \draw [thick,decorate,decoration={brace,amplitude=3pt,raise=0pt, mirror}] (1.0, -0.2) -- (1.5, -0.2);
    \node [below] at (1.3, -0.3) {$h$};
  \end{tikzpicture}
  \caption{Uniformly spaced fine-grid points and coarse-grid points with coarsening factor $m$. The $T_i$ are the C-points and form the coarse-grid, while the small hashmarks $t_i$  are F-points. Together, the F- and C-points form the fine-grid. }
  \label{fig:grids}
\end{figure}

We will first consider the two-grid scheme, consisting only of a fine grid, $\fg$, and a coarse grid, $\cg$, with coarsening factor $m$ in time. The multigrid method requires a coarsening scheme in time, intergrid transfer operators, and a coarse grid equation, defined as follows. To coarsen in time, label every $m$th time-point in $\fg$ a C-point and every other point an F-point, then $\cg$ is the set of C-points in $\fg$ (see Figure \ref{fig:grids}). A C-point, along with the following $m-1$ F-points, is called a coarse interval. For grid transfer operations, MGRIT uses \emph{injection}. For restriction, injection maps the values of $\vu$ at the C-points in $\fg$ to the points in $\cg$, and for interpolation, it maps the points in $\cg$ to the C-points in $\fg$. Following interpolation from $\cg$ to $\fg$, the solution on $\fg$ is relaxed using F-relaxation, which evolves the state at each C-point to the following F-points in each coarse interval using $\Phi$. F-relaxation may be viewed as interpolation, in that it fills in the F-points between each C-point. Importantly, since the coarse intervals are disjoint, F-relaxation can be done in parallel. We will not consider FCF-relaxation \cite{MGRIT14} here.

Coarsening in time induces a new problem posed on $\cg$, with ideal space-time operator $A_*$, in equation \eqref{eqn:splitting}, having a block subdiagonal consisting of $-\Phi^m$, where here $\Phi^m(\vu_i)$ is understood to mean the fine-grid propagation of the solution across the coarse interval to the next C-point, including forcing term $\vf$, e.g. $\Phi^2(\vv_{i-2}) = \Phi(\Phi(\vv_{i-2}) + \vf_{i-1})$, $\Phi^3(\vv_{i-3}) = \Phi(\Phi(\Phi(\vv_{i-3}) + \vf_{i-2}) + \vf_{i-1})$, etc. Solving this ideal coarse-grid equation gives the exact solution for each C-point, however, this requires as much work as solving the original fine-grid problem. Introducing the coarse grid time-stepping operator $\Phi_c$ which is designed to approximate the action of $\Phi^m$, MGRIT uses the approximation $A_c$, such that
\begin{equation} \label{eqn:splitting}
  A_* = \begin{bmatrix}
    I        &         &        &          &   \\
    - \Phi^m & I       &        &          &   \\
    & -\Phi^m & I      &          &   \\
    &         & \ddots & \ddots   &   \\
    &         &        & - \Phi^m & I
  \end{bmatrix} \approx A_c =
  \begin{bmatrix}
    I        &         &        &          &   \\
    - \Phi_c & I       &        &          &   \\
    & -\Phi_c & I      &          &   \\
    &         & \ddots & \ddots   &   \\
    &         &        & - \Phi_c & I
  \end{bmatrix}.
\end{equation}
Now, the equation $A_c(\vu) = R\vf = \vf_c$ is less expensive to solve than the fine-grid equation. In the case that $\Phi$ is derived from a continuous-time problem, then $\Phi_c$ is usually derived from a rediscretization of the continuous problem over the new coarse time-grid $\cg$. Deriving coarse operators for general $\Phi$ is an open problem \cite{ben_schroder_review,hyperbolic_coarse_grid}, and motivates this paper.
  
The approximation of $A_*$ by $A_c$ may be interpreted as a splitting method. Let $\vec{\tau}(\vu) = A_c(\vu) - A_*(\vu)$, then $A_*(\vu) = A_c(\vu) - \tau(\vu) = \vf_c$, and we immediately get the well-known $\tau$-correction form of Full Approximation Scheme (FAS) multigrid \cite{achi_77}:
  \begin{equation}
    A_c(\vv^{k+1}) = \vf_c + \vec{\tau}(\vv^{k}),
    \label{eqn:FAS_coarse_grid}
  \end{equation}
where $\vv^k$ is the approximate coarse solution on $\cg$ after $k$ multigrid iterations, and $\vec{\tau}_i = \Phi^m(\vv^k_{i-1}) - \Phi_c(\vv^k_{i-1})$. One two grid MGRIT iteration involves computing $\vtau(\vv^k)$ on $\fg$, injecting $\vf$ and $\vtau$ to $\cg$, solving \eqref{eqn:FAS_coarse_grid} sequentially for $\vv^{k+1}$, then interpolating to $\fg$ and applying F-relaxation.
$\vtau$ takes the form of a forcing term on the coarse grid, and it steers the solution toward the fine-grid solution, as well as ensuring that the exact fine-grid solution is a fixed point of the iteration.
The two grid algorithm is detailed in algorithm \ref{alg:MGRIT}.

\newcommand{\MGRIT}[1]{$\mathrm{MGRIT}_{#1}$}
Finally, we get the multigrid algorithm by applying the two-grid algorithm recursively. If \MGRIT{2} is the two grid algorithm, then we get the three-grid algorithm by replacing the sequential solve of the coarse grid equation \eqref{eqn:FAS_coarse_grid} with another application of \MGRIT{2}. Recursing this process gives us the V-cycle \MGRIT{m_l} algorithm, with $m_l$ levels.

\begin{algorithm}
  \label{alg:MGRIT}
  \caption{MGRIT 2 grid V-cycle, \MGRIT{2}$(\vv, \vf, m)$}
  \begin{algorithmic}
    \FOR{each C-point, $i = m, 2m, 3m, \dots, n$}
    \STATE $\vec{\tau}_i \gets \Phi^m(\vv_{i-m}) - \Phi_c(\vv_{i-m})$
    \ENDFOR
    \STATE{restrict $\vec{\tau}$, $\vf$ to the coarse grid, and solve:}
    \FOR{$i = m, 2m, 3m, \dots, n$}
    \STATE{$\vv_i \gets \Phi_c(\vv_{i-m}) + \vf_i + \vec{\tau}_i$}
    \ENDFOR
    \STATE{refine, then F-relax with $\Phi_h$}
  \end{algorithmic}
\end{algorithm}

\subsection{Motivation: Chaotic Problems and MGRIT}\label{sec:intro_chaos}
Chaotic systems are globally stable, deterministic systems which demonstrate sensitive dependence on initial conditions and system parameters, and which have trajectories that never settle down to a steady state solution or a periodic orbit for almost all initial conditions. 

To study MGRIT for chaotic systems, we will use the Lorenz system as a model problem. The Lorenz system is a three dimensional system of ODEs which is widely studied as an archetypal example of a chaotic system, and is given by
\begin{equation}
  \begin{cases}
    x' & = \sigma (y-x)    \\
    y' & = x(\rho - z) - y \\
    z' & = xy - \beta z
  \end{cases}.
\end{equation}
For the classical values of parameters $\sigma = 28$, $\rho = 10$, and $\beta = 8/3$, the Lorenz system is chaotic, with greatest \emph{Lyapunov exponent} of $\lambda_0 \approx 0.9$ \cite{strogatz}. This can be understood to mean that two trajectories differing only infinitesimally in initial conditions will, almost surely, diverge exponentially from each other in time with average rate $\lambda_0$. A system with $n_s$ dimensions has $n_s$ Lyapunov exponents, which are characteristic of the qualitative behavior of the system, and every chaotic system has a greatest Lyapunov exponent which is greater than zero.


PinT simulations of chaotic systems such as Lorenz are difficult because of two main problems. The first is that errors committed by the coarse operator will grow exponentially in time on the coarse grid. The other is that coarsening in time can cause serious qualitative changes in the behavior of the system caused by changes to the Lyapunov exponents. Therefore, the challenge is to form a coarse grid equation that is both locally precise and also captures the global qualitative behavior of the system.

We can quantify this difficulty with the condition number of the initial value problem \eqref{eqn:ivp}. First, we define \emph{Lyapunov time}, $T_\lambda = \frac{\ln(10)}{\lambda_0}$, to be the time it takes for a perturbation to a trajectory to grow by a factor of 10 \cite{strogatz}. In another sense, $T_\lambda$ is the time it takes for our numerical simulation of the system to lose one digit of accuracy. This provides an estimate for the condition number, $\kappa = \mathcal{O}(10^{T_f/T_\lambda})$, where $T_f = nh$ is the length of the time-domain.

\begin{figure}
  \centering
  \includegraphics[width=\textwidth]{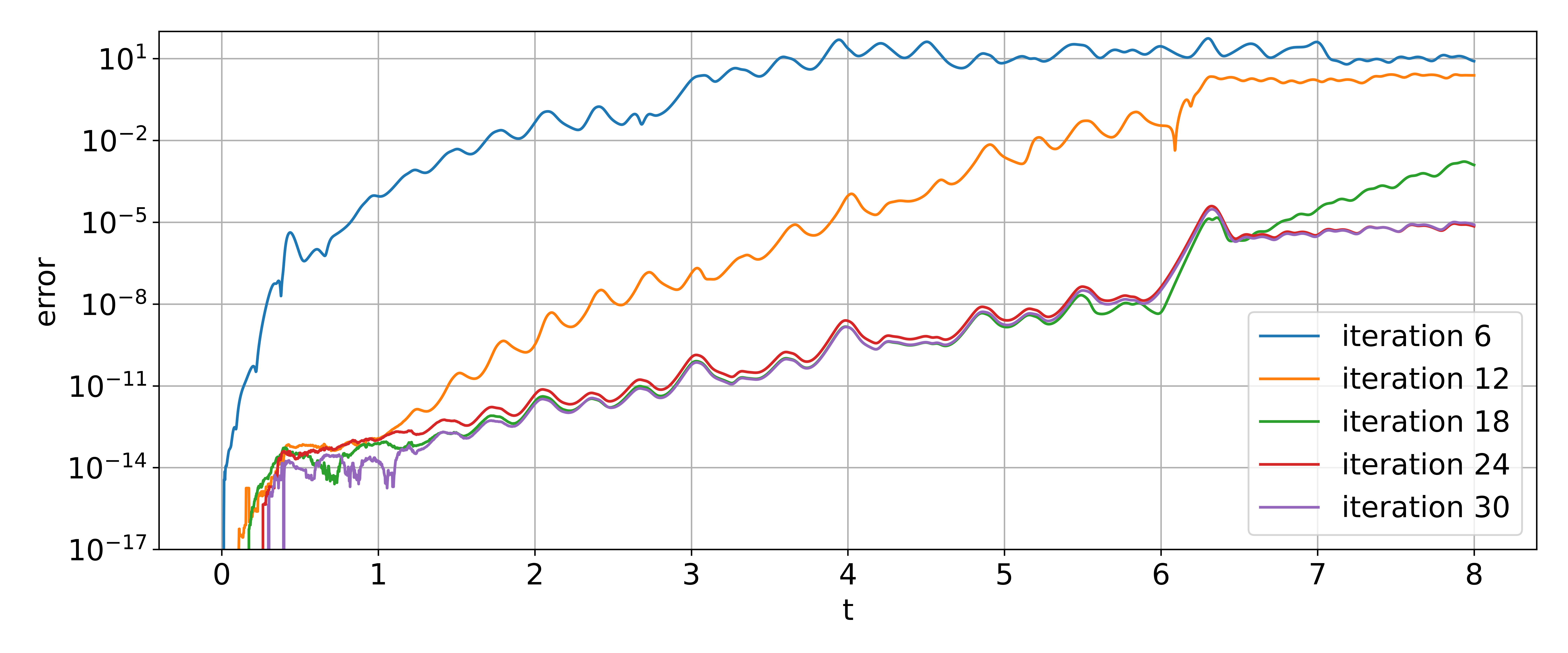}
  \caption{error over the time-domain (in units of Lyapunov time) for 30 MGRIT iterations on the Lorenz equation with $T_f = 8T_\lambda$. Convergence stalls after 24 iterations, due to the exponentially growing error.}
  \label{fig:MGRIT_stall}
\end{figure}

For example, Figure \ref{fig:MGRIT_stall} shows MGRIT performance for the Lorenz system with $T_f = 8T_\lambda$, giving us $\kappa = \mathcal{O}(10^8)$. The error, $\ve = \vv - \vu$ grows exponentially in time, with average rate $\lambda_0$. Thus, to converge to some tolerance, we would need the error to be made arbitrarily small at the beginning of the time-domain. In exact arithmetic, it can be shown that this is possible, but in practice we are limited by machine precision. If $\epsilon$ is the machine precision, then the maximum time domain size (in Lyapunov time) given some tolerance $Tol$ is approximately $T_f < \log_{10}(Tol/\epsilon)$, which suggests that for a tolerance of $10^{-10}$, we should only be able to converge on a time domain with $T_f < 6T_\lambda$.

\section{Main Contribution}\label{sec:main}

\subsection{\texorpdfstring{$\Delta$}{Delta} Correction}\label{sub:delta}
Let $\vu$ be a solution to the initial value problem \eqref{eqn:ivp}, for a system where the dynamics are chaotic. Then we linearize $\Phi$ at each point $\vu_i$ to get $F_i = D_u \Phi(\vu_i)$, where $D_u$ is the differential operator with respect to the spatial variables. The matrix $F_i$ is called the tangent linear propagator \cite{lyap_vecs}, as it describes the propagation of infinitesimal perturbations along $\vu$, i.e. $\Phi(\vu_i + \ve) - \vu_{i+1} \approx F_i \ve$. Thus, the propagation of a perturbation to the point $\vu_i$ is given by $W_n = F_n F_{n-1} \dots F_{i+1} F_i$. In the limit as $n \to \infty$, the time-average of the singular values of $W_n$ are equal to $\exp(\lambda_i)$, where the $\lambda_i$ are the Lyapunov exponents of the system. 

For systems where the underlying dynamics are chaotic, MGRIT is very sensitive to errors, no matter how small. This is because while the $\vec{\tau}$ correction drives the trajectory on the coarse grid toward that of the fine grid, the fine and coarse operators still have different tangent linear propagators, $F_i$ along that trajectory, meaning that even near convergence, there will be a significant mismatch between the Lyapunov spectrum on the coarse and fine grids. The $\Delta$ correction remedies this by using linearizations of the fine operator to update the coarse operator.
Let $\Phi$, $\Phi_c$, and $\Phi^m$ be defined as above, and $\vv = \vu - \ve$ be an approximate solution on the fine-grid. Now define
\begin{equation}
  \Delta_i := \left(D_u \Phi^m - D_u \Phi_c\right)(\vv_i)
  \label{eqn:delta}
\end{equation}
where $\Delta$ is a matrix valued function of $\vec{v}_i$, which encodes the difference between the linearizations of the ideal and coarse operators. Contrast this with the $\vec{\tau}$ correction, which encodes the difference between the \emph{values} of these two operators applied to $\vv_i$. The matrix $\Delta$ will naturally have the same number of dimensions as the number of spatial dimensions of the system.
We then use the computed $\Delta_i$ to form a correction to the time-stepper on each coarse interval:
\begin{equation}
  \Phi_\Delta(\vv_i) := \Phi_c(\vv_i) + \Delta_i \vv_i,
  \label{eqn:phi_coarse}
\end{equation}
which ensures that as $\vv$ approaches $\vu$, i.e. near MGRIT convergence, $D_u \Phi_\Delta$ approaches $D_u \Phi_*$.

Together with the $\vec{\tau}$ correction, which is computed at the same time, this gives the modified MGRIT algorithm \ref{alg:delta}, where the new additions are colored in red.

\begin{algorithm}
  \caption{MGRIT 2 grid V-cycle with $\Delta$ correction, $\Delta \mathrm{MGRIT}_2(\vv, \vf, m)$}
  \label{alg:delta}
  \begin{algorithmic}
    \FOR{each C-point, $i = m, 2m, 3m, \dots, n$}
    \STATE \textcolor{red}{$\Delta_i \gets D_u \Phi^m(\vv_{i-m}) - D_u \Phi_c(\vv_{i-m})$}
    \STATE $\vec{\tau}_i \gets \Phi^m(\vv_{i-m}) - \Phi_{\Delta_i}(\vv_{i-m})$
    \ENDFOR
    \STATE{restrict $\vec{\tau}$ and $\vf$ to the coarse grid, and solve:}
    \FOR{$i = m, 2m, 3m, \dots, n$}
    \STATE{$\vv_i \gets $\textcolor{red}{$\;\Phi_{\Delta_i}(\vv_{i-m})$} $+\;\vec{\tau}_i + \vf_i$}
    \ENDFOR
    \STATE{interpolate, then f-relax with $\Phi$}
  \end{algorithmic}
\end{algorithm}

Note that the first loop does not update the values of $\vec{v}$ at each time point, and may thus be done in parallel, while the loop on the coarse grid must be solved sequentially. Remember that, as before, the multigrid method replaces the forward solve on the coarse grid with a recursive call to the algorithm.

In the two-grid setting, $\Delta$-corrected MGRIT can be shown to be a generalization of Newton's method applied to the residual equation $\vr(\vv^k) = \vf - A_*(\vv^k) = 0$, with equivalence to Newton's method in the special case that $\Phi_c \equiv 0$. Thus, MGRIT with this $\Delta$ correction is expected to converge quadratically in certain regimes.

\subsection{\texorpdfstring{$\theta$}{Theta} method}
\begin{figure}
  \centering
  \includegraphics[width=0.8\textwidth]{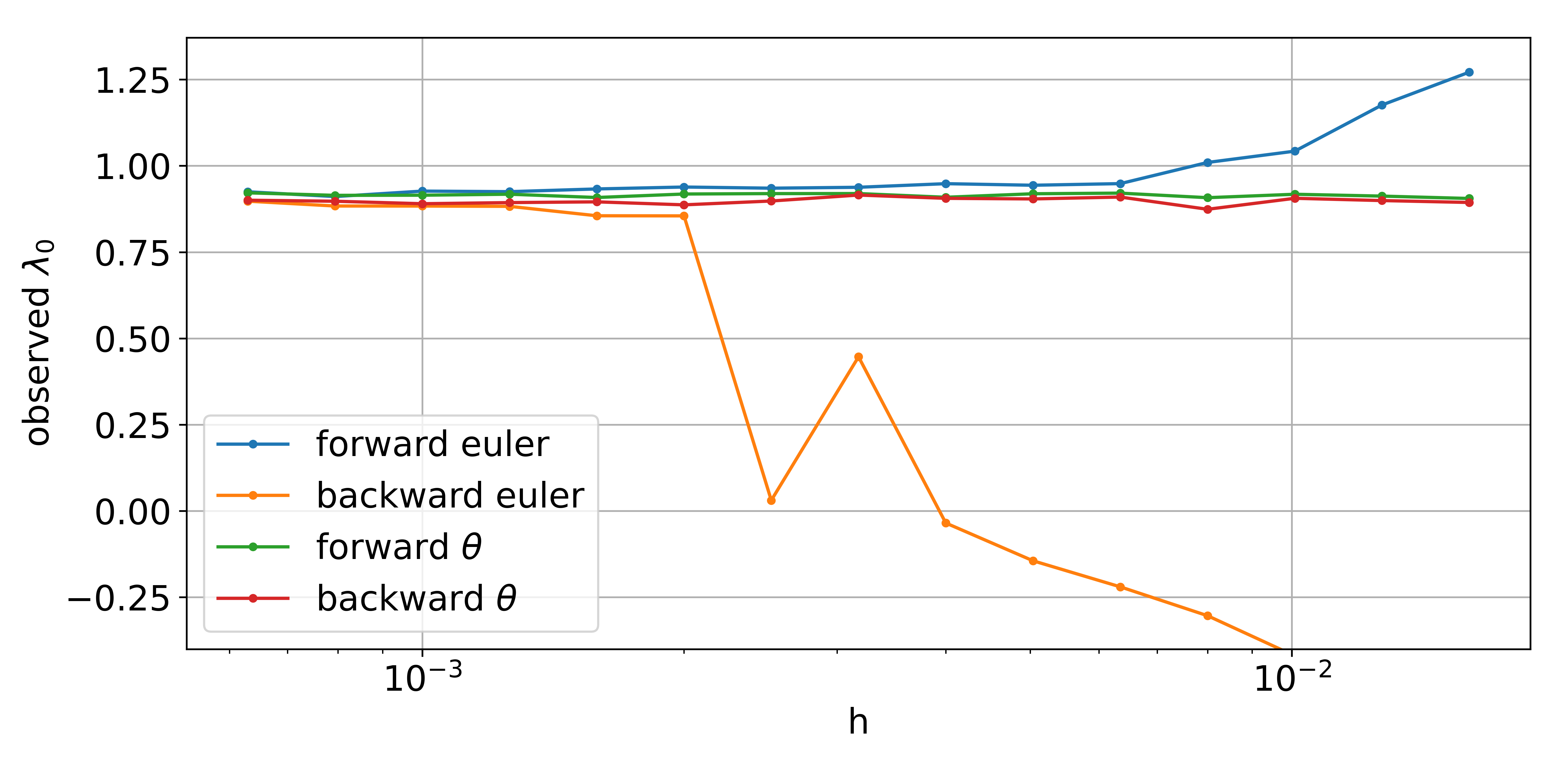}
  \caption{Plot of the observed greatest Lyapunov exponent $\lambda_0$ for different time-step sizes $h$, comparing forward Euler, backward Euler, as well as the $\theta$ method proposed here. For forward and backward Euler, coarsening in time changes the qualitative behavior of the system, while the $\theta$ method, using the ``forward" and ``backward" asymptotic values of $\theta$, preserves this sensitivity on coarse grids.}
  \label{fig:lyap_dependence}
\end{figure}

Another difficulty in solving chaotic systems with PinT is that coarsening in time can cause dramatic qualitative changes to the behavior of the system. For example, it is well documented that when using implicit Euler to solve the Lorenz equations, the measured greatest Lyapunov exponent decreases with increasing step size, $h$, meaning that for large $h$, a chaotic system can become artificially stablized. Conversely, using forward Euler, the Lyapunov exponent increases with increasing $h$, and the system appears \emph{more} chaotic on coarse grids \cite{what_good_are}. Figure \ref{fig:lyap_dependence} demonstrates this dependence for different time-stepping schemes applied to the Lorenz system, including for the $\theta$ method described here. Thus, if we seek a time-stepping scheme that preserves the qualitative behavior of the system on coarse grids, we should look for a scheme which lies somewhere between the implicit-explicit binary.

The $\theta$ method is a single-step time-stepping scheme designed to solve a discretized ODE given by $u_t = f(u)$. A single step is given by
\begin{equation}
  \vu_{i+1} = \vu_i + h [\theta f(\vu_i) + (1 - \theta)f(\vu_{i+1})],
  \label{eqn:theta}
\end{equation}
where the parameter $\theta \in [0, 1]$ gives control over the explicit/implicit character. For example, $\theta = 0$ corresponds to backward Euler, $\theta = 1$ to forward Euler, and $\theta = 1/2$ gives the second order trapezoid method. We use a simplified heuristic analysis to find $\theta$ values which work well in general.
Assume that the state variable $u$ is scalar, and that the fine-grid operator $\Phi$ is given by forward Euler. The coarse interval is then comprised of the points $u_0, u_1, \dots, u_m$, where $u_0$, $u_m$ are C-points, and $u_1, \dots, u_{m-1}$ are F-points, such that $u_{i+1} = u_i + hf(u_i)$. Compare the value of $u_m$ as computed by $\Phi^m$ and by $\Phi_\theta$, the theta method:
\[
  \Phi^m(u_0)       = u_0 + h\sum_{i=0}^{m-1} f(u_i)\;\;\;\text{and}\;\;\;\Phi_\theta(u_0)  = u_0 + mh[\theta f(u_0) + (1 - \theta)f(u_m)].
\]
We see that $m$ applications of forward Euler is equivalent to a left Reimann sum with $m$ points, which the $\theta$ method approximates with a weighted trapezoid rule. If we require that $f(u_m) \neq f(u_0)$, then we can solve for $\theta$ such that $\Phi_\theta(u_0) = \Phi^m(u_0)$, giving
\begin{equation}
  \theta = \frac{f(u_m) - \frac{1}{m}\sum_{i=0}^{m-1}f(u_i)}{f(u_m) - f(u_0)}.
  \label{eqn:theta_exact}
\end{equation}
This value of $\theta$ is computable, and gives us a coarse operator which is exact in the scalar case, however, it is not bounded, and can get large, especially at inflection points, where $u''(t) = 0$ and thus $f(u_n) \approx f(u_0)$. Further, computing $\theta$ this way requires several evaluations of the function $f$, which is likely not practical. We seek values of $\theta$ that work well in general. To simplify the notation, let $f_i$ denote $f(u_i)$. If we take the limit as $h \to 0$, we get the \emph{asymptotic values}, $\theta_m$:
\begin{equation}
  \theta_m  = 1 - \frac{1}{m}\sum_{i=1}^{m-1}\frac{ihf'_0 + \mathcal{O}(i^2h^2)}{mhf'_0 + \mathcal{O}(m^2 h^2)}
  \to  \frac{m+1}{2m}.
  \label{eqn:theta_ass}
\end{equation}
This formula for $\theta_m$ is then used to find a constant asymptotic value of $\theta$ for each coarse grid. Similarly, if we instead assume that the fine-grid operator is given by backward Euler, we get the asymptotic values $\theta_m = \frac{m - 1}{2m}$. Importantly, note that in either case, $\lim_{m \to \infty} \theta_m = \frac{1}{2}$ so $\theta_m \in (1/2, 1]$. Thus, this method should not be expected to be stable for arbitrarily coarse grids, and may be limited by the stability of the trapezoid rule, $\theta = 1/2$.
Surprisingly, although these asymptotic values of $\theta$ have been computed under the assumption that $u$ is scalar, they perform very well in the multivariate case as well. This can be seen in Figure \ref{fig:lyap_dependence} where the $\theta$ method using the forward and backward asymptotic values is compared to forward and backward Euler applied to the Lorenz system. We see that while the measured greatest Lyapunov exponent changes with increasing step size for forward and backward Euler, the $\theta$ method seems to preserve the Lyapunov exponent even on coarse time-grids.


\section{Numerical Results}\label{sec:num}

In the following experiments, we discretize the Lorenz system in time and solve using forward Euler's method on the fine grid. A coarsening factor of $m=2$ is used across all of the studies. When the $\theta$ method \eqref{eqn:theta} is used on the coarse grid, the values of $\theta$ are dependent on the grid level, $l$, and computed according to \eqref{eqn:theta_ass} with $m = 2^{l}$, where the fine-grid is assigned $l=0$, and the coarse grids are numbered $l=1, 2, \dots, m_l$. The implicit equation \eqref{eqn:theta} is solved numerically using Newton's method. When the $\theta$ method is not used on the coarse grid, forward Euler is used, with coarsened time-step size $m^{l}h$. First, we examine the convergence rates for the two grid algorithms on a small problem. Then we perform a refinement study and a time-domain size scaling study. Finally we explore the effect of adding more coarse levels for different problem sizes.

\begin{figure}
  \centering
  \includegraphics[width=0.9\textwidth]{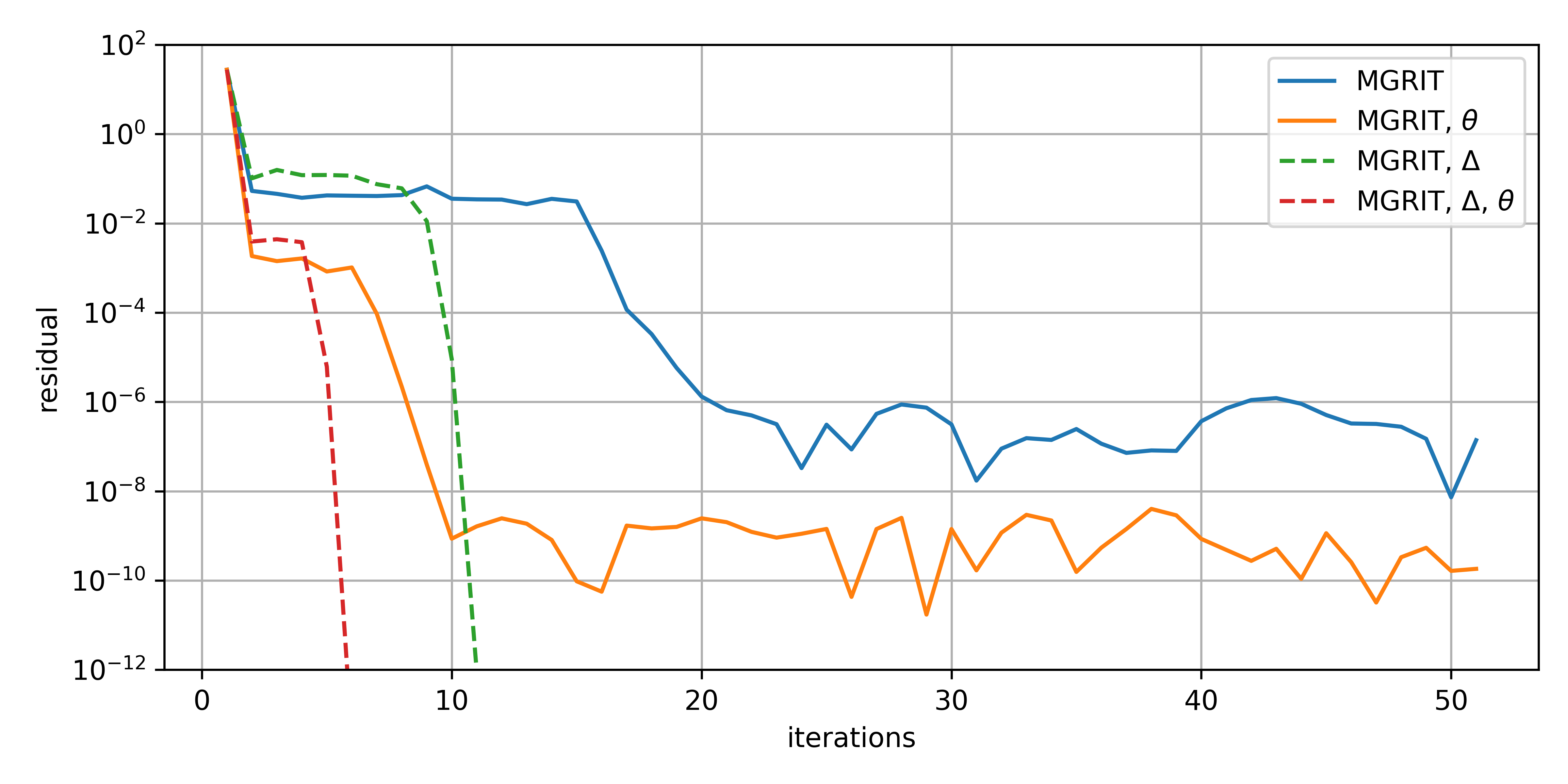}
  \caption{Residual history for each two grid algorithm applied to the Lorenz system with $T_f = 8T_\lambda$ and 8192 time-points. We see that \MGRIT{2} stalls on this time domain, as expected, but that we are able to converge using the $\Delta$ correction, and we see quadratic convergence for the $\Delta$ corrected algorithm.}
  \label{fig:convergence}
\end{figure}

Figure \ref{fig:convergence} plots the convergence history of the modified two grid MGRIT algorithms, solving the Lorenz system with $T_f = 8T_\lambda$. This experiment demonstrates that \MGRIT{2}, even using the $\theta$ method on coarse grids, stalls for long time-domains, which is expected given our pessimistic estimate on the loss of numerical precision from section \ref{sec:intro_chaos}. however, the $\Delta$ correction allows the method to converge.

In order to study the effect of varying time-step sizes $h$ on \MGRIT{2} performance, Table \ref{tab:refine_2} shows numbers of iterations required to reach a certain residual tolerance for the Lorenz system with fixed $T_f$ and increasing numbers of time-points $n_t$. While in all cases fewer iterations are needed for smaller $h$, compared to naive \MGRIT{2}, both the $\Delta$ correction and the $\theta$ method require roughly half as many iterations to converge, and when used together they require a quarter the number iterations. Further, we see that the $\theta$ method can account for instabilities on the coarse-grid which cause MGRIT to diverge.

\begin{table}
  \centering
  \caption{Iterations required for each two grid algorithm to converge to a residual tolerance of $10^{-10}$ for the Lorenz system with $T_f = 8T_\lambda$ and varying number of time-points $n_t$. `*' indicates that the algorithm diverged due to numerical instability}
  \begin{tabular}{lccccc}
    \toprule
    \multicolumn{1}{c}{}                   & \multicolumn{5}{c}{$T_f$, $n_t$}                                          \\
    \cmidrule(rl){2-6}
    Algorithm
    &4, 512&4, 1024&4, 2048&4, 4096&4, 8192\\
    \cmidrule(r){1-1} \cmidrule(rl){2-6}
    $\mathrm{MGRIT}_2$&*&44&22&15&12\\
    $\mathrm{MGRIT}_2$, $\theta$&19&13&9&7&6\\
    $\mathrm{MGRIT}_2$, $\Delta$&*&11&8&6&6\\
    $\mathrm{MGRIT}_2$, $\Delta$, $\theta$&8&6&5&4&4
    \\
    \bottomrule
  \end{tabular}
  \label{tab:refine_2}
\end{table}

Table \ref{tab:time_scaling_2} shows iteration counts for convergence of the two grid algorithm on the Lorenz system with increasing time-domain size $T_f$ and fixed time-step size $h$. For naive \MGRIT{2}, iteration counts increase linearly up until the critical time $T_f = 6T_\lambda$, after which naive \MGRIT{2} stalls. In contrast, the $\Delta$ correction and $\theta$ method greatly improve convergence for all time-domain sizes. Notably, the iteration counts for the $\Delta$ corrected algorithm are nearly flat, even for very long time-domain sizes, particularly when combined with the $\theta$ method coarse grid. 

\begin{table}
  \centering
  \caption{Iterations required for each two grid algorithm to converge to a residual tolerance of $10^{-10}$ for the Lorenz system with varying $T_f$ (in Lyapunov time) and varying number of time-points $n_t$ such that $h$ is constant. `-' indicates that the algorithm did not converge within 100 iterations}
  \begin{tabular}{lcccccc}
    \toprule
    \multicolumn{1}{c}{}                   & \multicolumn{6}{c}{$T_f$, $n_t$}                                  \\
    \cmidrule(rl){2-7}
    Algorithm
    &2, 4096&4, 8192&6, 12288&8, 16384&10, 20480&12, 24576\\
    \cmidrule(r){1-1} \cmidrule(rl){2-7}
    $\mathrm{MGRIT}_2$&10&13&17&64&-&-\\
    $\mathrm{MGRIT}_2$, $\theta$&4&5&6&7&41&-\\
    $\mathrm{MGRIT}_2$, $\Delta$&5&6&7&8&9&94\\
    $\mathrm{MGRIT}_2$, $\Delta$, $\theta$&3&4&4&5&5&48
    \\
    \bottomrule
  \end{tabular}
  \label{tab:time_scaling_2}
\end{table}

While \MGRIT{2} is not always used in practice, it is used as a stepping stone toward understanding the multilevel algorithm. Recall that \MGRIT{2} solves the coarse grid equation \eqref{eqn:FAS_coarse_grid} using a sequential solve, and \MGRIT{3} replaces this sequential solve with an application of \MGRIT{2} to the coarse grid, so we should expect that \MGRIT{3} will converge slower than \MGRIT{2}, and as we add more levels, this trend should continue. Thus we treat \MGRIT{2} as a best case scenario for the multilevel method.

Table \ref{tab:ml_scaling} demonstrates the effect that increasing the number of coarse grids has on MGRIT convergence. We see that adding a second coarse grid in \MGRIT{3} has a modest effect on convergence, while the jump from 2 coarse grids to 4, and 4 to 6 roughly double the iterations required for convergence, while quartering the problem size on the coarse grid.
The convergence of \MGRIT{7} with both the $\theta$ method and $\Delta$ correction is especially promising for $T_f=2$ and $T_f=4$, since the coarsest grid for those problem sizes is small (64 and 128 time-points respectively), and the iteration counts are similar to previous cases demonstrating parallel speedup \cite{MGRIT14}.

\begin{table}
  \centering
  \caption{Iterations required for each algorithm, using varying numbers of grids, to converge to a residual tolerance of $10^{-10}$ for the Lorenz system with varying $T_f$ Lyapunov time and varying number of time-points $n_t$ such that $T_f/n_t$ is constant. `-' indicates that the algorithm did not converge within 100 iterations, `*' indicates that the algorithm diverged due to numerical instability.}
  \begin{tabular}{lcccc}
    \toprule
    \multicolumn{1}{c}{}                   & \multicolumn{4}{c}{$T_f$, $n_t$}                                  \\
    \cmidrule(rl){2-5}
    Algorithm &2, 4096&4, 8192&6, 12288&8, 16384\\
    \cmidrule(r){1-1} \cmidrule(rl){2-5}
$\mathrm{MGRIT}_2$&10&13&17&64\\
$\mathrm{MGRIT}_3$&13&18&-&-\\
$\mathrm{MGRIT}_5$&26&-&-&-\\
$\mathrm{MGRIT}_7$&*&*&*&*\\
$\mathrm{MGRIT}_2$, $\theta$&4&5&6&7\\
$\mathrm{MGRIT}_3$, $\theta$&6&7&9&11\\
$\mathrm{MGRIT}_5$, $\theta$&10&13&19&63\\
$\mathrm{MGRIT}_7$, $\theta$&43&-&-&-\\
$\mathrm{MGRIT}_2$, $\Delta$&5&6&7&8\\
$\mathrm{MGRIT}_3$, $\Delta$&6&8&11&13\\
$\mathrm{MGRIT}_5$, $\Delta$&*&*&*&*\\
$\mathrm{MGRIT}_7$, $\Delta$&*&*&*&*\\
$\mathrm{MGRIT}_2$, $\Delta$, $\theta$&3&4&4&5\\
$\mathrm{MGRIT}_3$, $\Delta$, $\theta$&3&4&5&5\\
$\mathrm{MGRIT}_5$, $\Delta$, $\theta$&5&6&7&9\\
$\mathrm{MGRIT}_7$, $\Delta$, $\theta$&9&15&20&23\\
    \bottomrule
  \end{tabular}
  \label{tab:ml_scaling}
\end{table}

\section{Conclusions}\label{sec:conc}
Although simulating chaotic dynamical systems with PinT is inherently difficult due to the exponentially increasing condition number of the initial value problem, increased parallelism is necessary due to the importance of chaotic systems, as well as the increasing concurrency of high performance computing clusters. By adding a coarse grid correction based on the tangent linear propagator, and further improvement by exploiting the $\theta$ method on coarse grids, we observe nearly flat iteration counts for long time-domains in some cases. Our results demonstrate that PinT speedup may already be possible for the Lorenz system, which would be a first. This is supported by the fact that parallel speedups have been recorded for linear parabolic equations having similar problem sizes as those presented here. Additionally, while the techniques presented here are designed for chaotic problems, they will likely be effective for general nonlinear problems, which is future work.

\bibliographystyle{siamplain}
\bibliography{references_no_url}

\end{document}